\documentclass[11pt]{article} 
\usepackage[english]{babel}
 \usepackage{amssymb,amsbsy,amsmath,amsfonts,amscd}
\usepackage{graphics,graphicx} 
\usepackage{color} 
\usepackage{graphicx}
\def\R{{\mathbb R}} 
  
\def\v3{\vskip0.3cm \noindent}

 \newtheorem{rem}{\bf Remark}

\newtheorem{thm}{\bf Theorem\/} 
 
\newtheorem{lem}{\bf Lemma\/} 

\begin{document} 
\date{today} 
\centerline{\Large {\bf Semi-linear cooperative  elliptic systems }} 
\centerline{\Large {\bf involving Schr\" odinger operators: }} 
\centerline{\Large {\bf Groundstate positivity or negativity.}}
\vskip0.2cm 
\centerline{\bf B.Alziary - J.Fleckinger}

\vskip0.5cm \noindent 

{\bf Classification} 35J61, 35J10  
\vskip0.5cm \noindent
{\bf                     Abstract     }                      %
\vskip0.2cm \noindent
We study here the behavior  of the solutions to a $2\times 2$  semi-linear cooperative  system involving Schr\" odinger operators (considered in its variational form):
$$LU:=(-\Delta + q(x))U = AU+\mu U + F(x,U)\; \mbox{ in }\; \R^N$$
$$U(x)_{|x|\rightarrow \infty} \rightarrow 0 $$
where $q$ is a continuous  positive potential tending to $+\infty$ at infinity; $\mu$ is a real parameter varying near the principal eigenvalue of the system; $U$ is a column vector with components $u_1$ and $u_2$ and $A$ is a 
square cooperative  matrix with constant  coefficient.
$F$ is a column vector with components $f_1$ and $f_2$ depending eventually on $U$.
\section{Introduction}
We study here the behaviour of the solutions to a $2\times 2$  semi-linear cooperative  system involving Schr\" odinger operators (considered in its variational form):
$$LU:=(-\Delta + q(x))U = AU+\mu U + F(x,U)\; \mbox{ in } \; \R^N$$
$$U(x)_{|x|\rightarrow \infty} \rightarrow 0 $$
where $q$ is a continuous  positive potential tending to $+\infty$ at infinity; $U$ is a column vector with components $u_1$ and $u_2$ and $A$ is a 
square matrix with constant  coefficients; moreover $A$ is a cooperative matrix (which means that  its coefficients outside the diagonal are non negative).
$F$ is a column vector with components $f_1$ and $f_2$ depending eventually on $U$. The real parameter $\mu$ varies near the principal eigenvalue of the system and plays a key role. According to its position it determines not only the sign of the solutions but also their position w.r.t. the groundstate. 
\par \noindent
Such systems have been intensively studied (very often for $\mu = 0$) and mainly for Dirichlet problems defined on bounded domains ( \cite{FiMi1986}, \cite{FiMi1988}, \cite{FiMi1990},\cite{He1990},\cite{FHeT1995}, \cite{Sw1989},\cite{Am2005}, {\cite{AF2016}). 
When the whole  $\R^N$ is considered, as here, 2 cases are generally studied: either "Schr\" odinger systems"
(\cite{AbF1992},\cite{ABe2005},\cite{ACF1997}, \cite{AFTa1999}), that is system involving    Schr\" odinger operators, as here, or systems with a weight tending to $0$ (\cite{FSe1995},\cite{AFLeWe2012}).
 It is also possible to consider a combination of these 2 problems with a potential $q$ and a weight $g$ :$$LU:=(-\Delta + q(x))U = g(x) AU+\mu  g(x) U + F(x,U)\;
 \mbox{ in } \; \R^N$$ as far as $\dfrac{g}{q} $ tends to $0$ at infinity which is the condition for having some compactness and therefore a discrete spectrum.
\par \noindent
The first results on  Schr\" odinger  systems, when $F$ does not depend on $U$ (linear systems)  deal with cooperative systems and with the Maximum Principle ({\bf MP}) that is:
\par \noindent  {\it "If the data $F$ is non negative, $\neq 0$, then, any solution $U$ is non negative"}. 
\par \noindent 
As for the case of one equation, this Maximum Principle holds for a parameter $\mu< \Lambda^*$, where $\Lambda^*$ is the principal eigenvalue of the system, which means that 
$LU-AU-\Lambda^* U=0$ has  a non zero  solution which does not change sign. 
\par \noindent
For the classical case of  an equation defined on a bounded domain with zero boundary conditions, $-\Delta u=\mu u + f(x), \; f>0$ , Cl\'ement and Peletier  \cite{ClPe1979} 
 have shown that the solution
 $u$ changes sign as soon as $\mu$ goes over $\lambda_1$, the first eigenvalue of the Dirichlet Laplacian defined on $\Omega$. More precisely there exists a small  positive $\delta$,
 depending on $f$, such that for all 
$ \mu \in (\lambda_1, \lambda_1 + \delta), \; u<0.$ This phenomenon is known as "Anti-maximum Principle" ({\bf AMP}).  
\par \noindent
In our present case, where we have no boundary,  we have improved these results  giving not only the  sign of the solutions but  also comparing  the solutions with the groundstate 
(principal eigenfunction); it is   what we call  "groundstate positivity"({\bf GSP})  (resp. negativity) (resp. {\bf GSN}). 
We extend in particular previous results established in  \cite{AF2017} for linear systems   to some semi-linear cooperative systems. For being not excessively  technical,
 we limit our study  to radial potentials and cooperative systems. Extensions to more general cases 
will appear somewhere else. 
\vskip0.3cm \noindent
Our paper is organized as follows:
\par \noindent
We  recall first some previous results of the linear case  that we use.
Then we study  a semi-linear equation.
Finally we study a cooperative 
 semi-linear system.

\par \noindent


\section{Linear Case: one equation}
We shortly recall the case of a linear equation with a parameter $\mu$
varying near  the principal eigenvalue of the operator.
 

$$ Lu:= (-\Delta +q(x)) u \, = \, \mu u + f(x) \mbox{ in }\R^N, \leqno(E)$$ 
$$ \lim_{|x| \rightarrow + \infty}  u(x) \, = \, 0.$$
$$
(H_q) \quad  q \mbox{ is a positive continuous potential tending to $+\infty$ at infinity}. \quad \quad
\quad \quad\quad \quad\quad \quad\quad \quad\quad \quad
$$
We seek $u$ in $V$ where
$$V:=\{ u \in L^2(\R^N)\, s.t. \; \|u\|_V = \big(\int|\nabla u|^2 + q(x)u^2 \big)^{1/2} <\infty\}.$$
If $(H_q)$ is satisfied, the embedding of $V$ into $L^2(\R^N)$ is compact (see $e.g.$ \cite{F1981},\cite{Ed-Ev}). 
Hence $L$ possesses an infinity of 
eigenvalues tending to $+\infty$: 
$$0< \lambda_1 < \lambda_2 \leq ....\leq \lambda_k \leq ... \, , \; \lambda_k \rightarrow +\infty \; {\rm  as} \; k\rightarrow \infty.$$
\par \noindent{\bf Notation : $(\Lambda, \phi)$} $\,$ We set  from now on
 $\Lambda:= \lambda_1 $ the smallest one (which is positive and simple) and $\phi$ 
the associated eigenfunction, positive and with $L^2$-norm $\|\phi\|=1$.
\par \noindent 
It is classical (see $e.g.$ \cite{ReSi})  that if $f \geq 0, \neq 0, $ and $\mu <\Lambda$, 
there exists exactly one solution which is positive: the positivity is "improved",
 or in other words, the (strong) maximum principle {\bf  (MP)} is satisfied:
$$f\geq 0,\, 
 \not \equiv 0\; \Rightarrow\;  u>0. \leqno(MP)$$ 
Lately, as said above,    another notion   has been
defined (\cite{AFTa2001},\cite{ATa2007}, \cite{Le2010}) 
the "groundstate positivity" ({\bf GSP}) 
(resp. " negativity"  ({\bf GSN})) which means that,  there exists $k>0$   such that the solution 
$u>k\phi$ (GSP)  (resp. $u<-k\phi$ (GSN)). 
\par \noindent
We also say shortly  "fundamental positivity" or" negativity", or also "$\phi$-positivity" or "negativity". 
Indeed these properties are more precise than MP or AMP. But for proving them, it is necessary to have
 a potential growing fast enough, a potential  with a super quadratic growth.
\par \noindent
In  \cite{ATa2007} a class   ${\cal P}$  of radial potentials is defined:
\begin{equation} \label{P}
{\cal P}:=\{Q \in {\cal C}(\R_+,\R^*_+) /  \exists R_0>0, Q'>0 \, a.e. \, \mbox{on} \, [R_0, \infty), \,  
\int_{R_0}^{\infty} Q(r)^{-1/2} < \infty\}.\end{equation}
The last inequality  holds precisely  if $Q$ is growing sufficiently fast, indeed faster than  $r^2$ (the harmonic oscillator). 
In this paper we consider only  a
radial  potential $q \in {\cal P}$. Note that our proof is valid for more general potentials, in particular for perturbations of radial potential \cite{AFTa2007} or \cite{ATa2007} . We assume here 
$$
(H'_q) \quad   q \mbox{ is radial and  is in }{\cal  P}\quad\quad\quad \quad\quad\quad \quad\quad\quad \quad\quad\quad \quad\quad\quad \quad\quad\quad\quad\quad \quad\quad\quad \quad\quad\quad \quad\quad\quad \quad\quad\quad \quad\quad\quad\quad\quad \quad\quad\quad \quad\quad\quad \quad\quad\quad \quad\quad\quad \quad\quad
$$

\begin{rem}:  Note that since $q$  is in $ {\cal P}$ it satisfies $(H_q)$. 
\end{rem} 
\par \noindent
On $f$ we assume 
$$
(H^*_f) \quad  f \in L^2(\R^N), \quad f^1=\int f \phi >0. \quad \quad \quad \quad \quad \quad\quad \quad \quad\quad \quad \quad\quad \quad \quad\quad \quad \quad\quad \quad \quad\quad \quad \quad
$$
For having more precise estimates on $u$,  in particular the "groundstate negativity" {\bf (GSN)} , 
we have to define another   set $X$ in which $f$ varies,  the set of 
"groundstate bounded functions":
\begin{equation} \label{X}
X:=\{ h \in L^2(\R^N): \, |h|/\phi \in L^{\infty}(\R^N)\}, \end{equation}
equipped with the norm $\|h\|_X=ess\sup_{\R^n}(|h|/\phi)$.
\begin{thm}\label{LIM}: 
Assume $(H'_q)$ and $(H^*_f)$, $f\in X$. 
For  $  \mu < \Lambda$ or $\Lambda<\mu  < \lambda_2$ 
there exists $\delta>0$ (defined below)  depending on $f$  and a positive constant $C$, depending on $f$  such that
if $0<|\Lambda - \mu| < \delta$, 
$$\Lambda - \delta < \mu < \Lambda \; \Rightarrow \;  u \geq \frac{C}{\Lambda - \mu} \phi >0, $$
$$ \Lambda < \mu < \Lambda + \delta \; \Rightarrow \;  u \leq \frac{C}{\Lambda - \mu} \phi <0.$$
\end{thm}
{\bf   Proof of Theorem \ref{LIM}: }$\;$
Decompose now $u$ and $f$ in  $(E)$ on $\phi$ and its orthogonal:
$$u= u^1 \phi +  u^{\perp}\, ; \; f=f^1 \phi + f^{\perp}; \;  u^1=\int u \phi, \; \int  u^{\perp} \phi = \int  f^{\perp} \phi_,= 0;$$
we derive from Equation $(E)$
\begin{equation}\label{equ} (L-\mu) u^1\phi =   (\Lambda - \mu) u^1 \phi = f^1\phi\, , 
\;    L u^{\perp} = \mu u^{\perp} + f^{\perp} .\end{equation}
Choose $ \mu < \Lambda$ or $\Lambda<\mu  < \lambda_2$ .
From the first equation we derive    
 $$u^1= \frac{f^1}{(\Lambda - \mu)} \, \rightarrow \, \pm \infty \, \mbox{ as }
\, (\Lambda - \mu)  \rightarrow 0.$$
 By use of Theorem 3.2 (c)  in \cite{AFTa2007} or \cite{ATa2007}, 
 we know that the restriction of the resolvent $(L-\mu)^{-1} $ to $X$ is bounded from $X$ into itself. The following lemma is a direct consequence of this result as it is shown in the proof of the Theorem 3.4 in \cite{AFTa2007}.

 \begin{lem}\label{delta}:
 There exists  $\delta_0$ small enough and  
 there exists a constant $c_0$ (depending on $\delta_0$) such that for all $\mu$ with 
 $\Lambda - \delta_0 < \mu < \Lambda$ or $\Lambda < \mu < \Lambda + \delta_0 < \lambda_2$,  
 $$-c_0 \|f^{\perp}\|_X  \leq    \|u^{\perp}\|_X \leq c_0 \|f^{\perp}\|_X. $$
\end{lem}
Finally we take in account  Lemma \ref{delta} and (\ref{equ}): 
$$\|u^{\perp}\|_X \leq c_0 \|f^{\perp}\|_X \;\;  {\rm and } \; \, 
u=\frac{f^1}{\Lambda - \mu} \phi + u^{\perp};$$
for $|\Lambda - \mu| \rightarrow 0, $  $\frac{f^1}{\Lambda -\mu} \phi \rightarrow \pm \infty$ when $u^{\perp}$ stays bounded.  
Hence, for $|\Lambda - \mu|$ small enough, more precisely for $|\Lambda - \mu| < \delta_1(f):= \frac{f^1}{c_0 \| f^{\perp}\|_X}$, we have 
$$\frac{f^1}{|\Lambda - \mu|} > c_0  \|f^{\perp}\|_X.  $$ 
We deduce that  Theorem 1 is valid for  $\delta:= \min\{\delta_0, \delta_1(f)\}$.

\section{Semi-linear  Schr\" odinger equation}
We study now the case of a semi-linear  equation. We first  obtain bounds for the  solutions, if they exist and then
we show their existence via the method of "sub-super solutions". 
Finally, with additional assumptions, we prove the uniqueness of them. 
\par \noindent
Consider  the semi-linear Schr\" odinger equation (SLSE)  
$$ Lu:= (-\Delta +q(x)) u \, = \, \mu u + f(x,u) \mbox{ in } \R^N, \leqno(SLSE)$$ 
$$ \lim_{|x| \rightarrow + \infty}  u(x) \, = \, 0.$$
We assume that the potential $q$ satisfies $(H'_q)$ and we denote as above by $(\Lambda, \phi)$ the principal eigenpair with $ \phi >0$.
\par \noindent
We work in $L^2(\R^N)$ and we consider the problem in its variational formulation.
We seek $u$ in $V$ for a  suitable   $f$. 
\par \noindent
We assume that $f$ satisfies :
\par \noindent
$(H_f)$ $f:\R^N \times \R \rightarrow \R$  is a Caratheodory function $i.e.$ the function $f(\bullet,u)$ is Lebesgue measurable in 
$\R^N$, for every $u(x)\in \R$ and the function $f(x,\bullet)$ is continuous in $\R$ for almost every $x\in \R^N$. 
Moreover, $f$ is such that   
$$\forall u \in L^2(\R^N ), \;\;  f(.,u) \in L^2(\R^N), \leqno\quad (i)$$
$$ \exists \, \kappa>0  \quad s.t.  \quad   \forall  u \, \in\, V,  \;\; f(x,u) \, \geq \, \kappa \phi(x) >0  \leqno \quad (ii)$$
$$  \exists K > \kappa >0 \quad s.t.  \quad  \forall  u \, \in\, V, \;\;  f(x,u) \, \leq \, K \phi(x). \leqno \quad (iii)$$

\par \noindent
Later  we also suppose 
$$\forall x \in \R^N, \; u \rightarrow \frac{f(x,u)}{|u|} \, {\rm is\,  strictly \, decreasing } \;   \leqno(H'_f)$$

\par \noindent
\begin{rem}:
Note that, by $(ii)$ and $(iii)$, for any $u\in V$,  $f(.,u) \in X$ and hence the solutions, if they exist, are in $X$. 
\end{rem}

\par \noindent
Let a parameter  $\mu $  be given, with $|\mu - \Lambda|$ ``small enough''. In this section we prove groundstate positivity and negativity for the semi-linear Schr\" odinger equation.
\begin{thm}\label{T1}:  If $(H'_q)$
and  $(H_f)$ are satisfied , then there exists $\delta(f) >0$
 ($\delta=\delta(f):= \min\{\delta_0, \delta_1'(f):= \frac{\kappa}{c_0 K}\}$  where $\delta_0$ and $c_0$  are given in Lemma 1)
 such that, for $0< |\mu - \Lambda| < \delta$ there exists a solution $u$ to $(SLES)$ such that
 $$\|u\|_X \leq \frac{K}{|\Lambda - \mu| } +2c_0 K .$$
 Also 
\par 
 - for $\Lambda - \delta < \mu < \Lambda$, $u> \frac{\kappa}{\Lambda - \mu} \phi>0$, 
\par
- for $\Lambda < \mu < \Lambda + \delta <\lambda_2$,   $u < \frac{K}{\Lambda - \mu}\phi<0$. 
\par \noindent
Moreover if $(H'_f)$ is satisfied, the solution to  $(SLSE)$ is unique.
\end{thm} 
\begin{rem} \label{RT1}
If $(ii)$ does not hold, for $\mu < \Lambda$, there exists a solution $u$ such that 
$$\|u\|_X \leq \frac{K}{|\Lambda - \mu| } +2c_0 K .$$ The existence is classical
($e.g.$ \cite {ACF1997}) and the estimate follows from the proof below.\end{rem}
\vskip0.2cm \noindent
{\bf Proof of Theorem \ref{T1}:}$\;$
\par \noindent 
We do the proof in $3$ steps:  first  maximun and anti-maximum principles, secondly  existence of the solution such that  $u> \frac{\kappa}{\Lambda - \mu} \phi>0$ for $\Lambda - \delta < \mu < \Lambda$ and such that $u < \frac{K}{\Lambda - \mu}\phi<0$, for $\Lambda < \mu < \Lambda + \delta$, and thirdly the uniqueness.

\vskip0.2cm \noindent
{\bf  Step 1. Maximun and anti-maximum principles }
\par \noindent
We prove the  positivity or negativity of the solutions exactly as for the linear case,  but, since $f$ depends on $u$ we have to show that  
$\delta$ (which depends on $f$ in the linear case) is now uniform. This follows from hypotheses 
$(ii)$ and $(iii)$.  
\par \noindent
Let $u$ be a solution to $Lu=\mu u + f(x,u)$. For this $u$, 
set $$f^1(u)=\int f(x,u) \phi(x) dx\, , \;  f^{\perp}(x, u) = f(x,u) - f^1(u) \phi(x).$$ 
 Also $u^1=\int u \phi(x) dx$ and 
$u^{\perp} = u - u^1 \phi$. 
\par \noindent 
Note that, always by $(ii)$ and $(iii)$, $0< \kappa \leq f^1(u) \leq K$. 
\par \noindent
With this decomposition, reporting in $(SLSE)$, we obtain 2 equations: 
$$ (L-\mu) u^1\phi =   (\Lambda - \mu) u^1 \phi = f^1 \phi\, , 
\;\,     L u^{\perp} = \mu u^{\perp} + f^{\perp} .$$ 
Choose $ \mu < \Lambda$ or $\Lambda<\mu  < \lambda_2$ . 
From the first equation we derive  
 $$u^1= \frac{f^1}{(\Lambda - \mu)} \, \rightarrow \, \pm \infty \, as
\, (\Lambda - \mu)  \rightarrow 0.$$
Now we proceed exactly as for the linear case.
 By use of Theorem 3.2 (c)  in \cite{AFTa2007} or \cite{ATa2007}, 
 we know that the restriction of the resolvent $(L-\mu)^{-1} $ to $X$ is bounded from $X$ into itself. 
 So by $(iii)$ and by Lemma \ref{delta} 
there exists a $\delta_0$ small enough and  
 there exists a constant $c_0$ (depending on $\delta_0$) such that for all $\mu$ with $|\Lambda - \mu| < \delta_0$,  
 $$\|u^{\perp}\|_X \leq c_0 \|f^{\perp}(x,u)\|_X\leq c_0\|f(x,u)-f^1(u) \phi(x)\|_X\leq 2c_0K. $$
 Write now
 $$u=\frac{f^1(u)}{\Lambda - \mu} \phi + u^{\perp}$$
 Hence $\|u\|_X \leq \frac{f^1(u)}{|\Lambda - \mu|} + \|u^{\perp} \|_X  \leq \frac{K}{|\Lambda - \mu|} + 2c_0 K $.
 For $|\Lambda - \mu| \rightarrow 0, $  $\frac{f^1}{\Lambda - \mu} \phi \rightarrow \pm \infty$ when $u^{\perp}$ stays bounded.  
For $|\Lambda - \mu|$ small enough, that is here
 $|\Lambda - \mu| < \delta_1'(f):= \frac{\kappa}{2 c_0 K}$, we get  (since $f^1>0$) 
 $$\frac{f^1}{|\Lambda - \mu|}\geq \frac{\kappa}{|\Lambda - \mu|} >2c_0K \geq  c_0  \|f^{\perp}\|_X.  $$
 Finally  Maximum and anti-maximum principles are valid for \\ $\delta(f) := \min\{\delta_0, \delta_1'(f)\}$.
 

\vskip0.2cm \noindent
{\bf Step  2. Existence of solutions
}
\par \noindent
We prove the  existence of solutions by Schauder fixed point theory; for this purpose we need some classical elements: a set ${\cal K}^{\pm}$  constructed with the help of sub-super solutions and a compact operator $T$  acting in    ${\cal K}^{\pm}$ such that 
${\cal K}^{\pm}$ stays invariant by $T$: $T({\cal K}^{\pm}) \subset {\cal K}^{\pm}$.
\vskip0.2cm \noindent
\underline{ 1:  "Sub-super solution" :} 
\par 
$\bullet \;$ Case $\Lambda - \delta < \mu < \Lambda$. 
\par \noindent
Obviously, by $(ii)$,    $u_0 = \frac{\kappa}{\Lambda - \mu}   \phi>0$  is a subsolution:
$$L(u-u_0)= \mu (u-u_0)+ f - (\Lambda - \mu) u_0 =\mu(u-u_0) +f-\kappa \phi$$ and by $(ii)$ and 
GSP, $u-u_0\geq 0$.
\par \noindent
Analogously  $u^0=  \frac{K}{\Lambda - \mu}   \phi \, >0$
( $K$  given in $(iii)$)  is a supersolution : 
 $$Lu^0 = \frac{\Lambda}{\Lambda - \mu} K\phi = \Lambda u^0= \mu u^0 + (\Lambda - \mu) u^0.$$
\begin{rem}:  The  sub- and supersolutions tend to $+\infty$ as $\mu \nearrow \Lambda$. \end{rem}

\par 
$\bullet \;$ Case $\Lambda < \mu < \Lambda +  \delta < \lambda_2 $.
 $v^0=  \frac{\kappa}{\Lambda - \mu}   \phi \, <0$ is a supersolution. Indeed
$$L(v^0- u) \, = \mu (v^0-u) +  \kappa \phi -f$$
 and by $(H_f)$ and  the anti-maximum $0> v^0\geq u$. 
\par \noindent 
Analogously, $v_0= \frac{K}{\Lambda - \mu}   \phi \, <0$ is a subsolution. 
\begin{rem}:  The sub- and supersolutions tend to $-\infty$ as $\mu \searrow \Lambda$. \end{rem}
\par \noindent
\begin{rem}:  Obviously,  $u_0 < u^0$ for $\Lambda -\delta < \mu < \Lambda $ (resp. $v_0<v^0$ for $\Lambda < \mu < \Lambda +  \delta $). \end{rem}
\underline{  2: The operator $T$}
\par \noindent
We define $T: u\in L^2 \longrightarrow w=Tu\in V, $ where $w\in X$ is the unique solution to
$ Lw= \mu w + f(x,u)$.
\vskip0.2cm  \noindent
\underline{ 3: $\;$  The  invariant set ${\cal K^+}:=[u_0, u^0 ]$} for $\Lambda -\delta < \mu < \Lambda $ (resp. ${\cal K^-}:=[v_0, v^0 ]$ for $\Lambda < \mu < \Lambda +  \delta $).
\par \noindent
If  $\mu < \Lambda$, by  the maximum principle and the hypothesis $(iii)$ ,  
$u\leq u^0 $ implies $ w\leq u^0$. Indeed,
$$L(u^0 -w) =   \mu (u^0 -w) +(\Lambda - \mu) u^0 -f(x,u) =   \mu (u^0 -w) +  K\phi - f(x,u); $$
since, by $(iii)$,   $ K\phi - f(x,u)\geq 0$, we apply the maximum principle and hence $w\leq  u^0$. 
The  3 other cases lead to analogous calculation. 
\vskip0.2cm \noindent
\underline{ 4: $T$ is compact in $X$.}
\par \noindent
First note that ${\cal K^+}\subset X$  (resp. ${\cal K^-}\subset X$). 
$Lw-\mu w = f(x,u)$ can also be written $ w= (L-\mu I)^{-1} f(x,u) = T(u)$. Since by \cite{ATa2007} ,\cite{AFTa2007},  the resolvent 
$ R(\mu):= (L-\mu I)^{-1}$ is compact  in $X$ for $\mu \in (\Lambda - \delta, \Lambda)$ or 
$(\Lambda, \Lambda + \delta)$,   and  since $F:u \rightarrow f(x,u)$ is continuous, $T= R(\mu)F$ is compact.
\vskip0.2cm \noindent
We deduce from Schauder fixed point theory that there exists a solution to $(SLSE)$ in $\cal K^+$, (resp. in $\cal K^-$ ).
\vskip0.2cm \noindent
{\bf Step 3. Uniqueness }$\;$
\par \noindent

For proving  uniqueness we follow  \cite{BrOs}, p.57. First we assume not only $(H_f)$ but also $(H'_f)$. 
Assume that $u$ and $v$ are two   solutions: 
\[Lu = \mu u + f(x,u) \;, \;\,Lv = \mu v + f(x,v) \]
The solutions are in $X$ and we have shown that  $ u,v>u_0>0$  for $\Lambda -\delta < \mu < \Lambda $ (resp. $u,v<v^0<0$ for $\Lambda < \mu < \Lambda +  \delta $).
 Hence we can write
\[ \frac{Lu}{u} = \mu + \frac{f(x,u)}{u}\,; \; \frac{Lv}{v} = \mu + \frac{f(x,v)}{v}.\]
By subtraction  $q(x)$ and  $\mu$ disappear. Multiply by $u^2-v^2$ and integrate. 
$$\int [\frac{-\Delta u}{u} + \frac{\Delta v}{v} ][u^2-v^2] \,=\, \int[\frac{f(x,u)}{u}-\frac{f(x,v)}{v}][u^2-v^2]; $$
the last term is non positive by $(H'_f)$.
\par \noindent
 We transform exactly as in \cite{BrOs} the first term.
\[
\int [\frac{-\Delta u}{u} + \frac{\Delta v}{v} ][u^2-v^2] \,=
\, \int |\nabla u - \frac{u}{v} \nabla v|^2 +  |\nabla v - \frac{v}{u} \nabla u|^2 =\]
\begin{equation}\label{brezis}
 \int|v\nabla\big(\frac{u}{v}\big)|^2+|u\nabla\big(\frac{v}{u}\big)|^2 \,\geq 0; 
\end{equation}  

therefore both terms are equal to $0$ and 
$$u^2-v^2=0  \, \Rightarrow \, u=v \,a.e.; $$ by regularity, $u=v$. 


\section{Semi-linear cooperative system}
\noindent
We extend here  to a class of semi-linear systems previous results shown in \cite{AF2017} where 
linear systems of the form $LU=\mu U +AU +F(x)$ are studied. 
\par \noindent
We study for $a>0$, $b>0$, $c>0$ 
\[ \left \{ \begin{array}{lcl}
Lu_1&=&(\mu +   a) u_1 +bu_2 + f_1(x,u_1)\\
 Lu_2&=& cu_1 + (\mu +d) u_2 +f_2(x,u_2)\end{array} \right .\, \; in \; \R^N ,.\leqno(S)\]
$$ \; u_1(x), u_2(x)_{|x|\rightarrow \infty} \rightarrow 0 .$$
We write shortly $LU= \mu U+AU+F(x,U),$ where $A$ is the cooperative matrix with components $a,b,c,d$:  $$A  =\left(\begin{array}{cccc}
a&b\\c &d
 \end{array}
 \right).$$ 
\par \noindent
{\bf Notation  $(\xi_1, Y)$:}  
Denote $\xi_1$ the  largest eigenvalue of $A$ (the other one being  denoted by $\xi_2$);    $Y$  is the  eigenvector associated with $\xi_1$: 
$$AY=\xi_1 Y.$$
$$\xi_1= \frac {a+d+\sqrt{(a-d)^2+4bc}}{2}.$$
An easy calculation shows that  $( L-A) (Y\phi) =(\Lambda - \xi_1) Y \phi ;$
moreover here $Y\phi$ is with components which do not change sign: we choose   both components of $Y$   positive:
$$y_1 = b>0 \, , \; y_2 = \frac{d-a + \sqrt{(a-d)^2+4bc}}{2}>0.$$ 

\par \noindent
{\bf Notation $\Lambda^*$:} 
 $\Lambda^*:=\Lambda - \xi_1$  is the principal eigenvalue of System $(S)$ with associated eigenvector $Y \phi$:
 $$( L-A) (Y\phi) =(\Lambda - \xi_1) Y \phi =\Lambda^* Y\phi .$$
\par \noindent
{\bf Hypotheses:}
We assume 
\par \noindent
 $(H_A)$ $\,$   $A$ is a $2 \times 2$  cooperative matrix with positive coefficients outside the diagonal. 
\par \noindent
$(H_{F}):$  $\;$  $f_1, f_2 : \R^N \times \R \rightarrow \R$ are  Caratheodory function $i.e.$ the functions
 $f_1(\bullet,u_1)$ or $f_2(\bullet,u_2)$  are Lebesgue measurable in 
$\R^N$, for every $u_1(x)$ or $u_2(x)$  in $\R$ and the functions $f_1(x,\bullet)$, $f_2(x,\bullet)$  are continuous in
 $\R$ for almost every $x\in \R^N$. 
Moreover, $f_1$, $f_2$ are such that   
$$\forall u_1, u_2  \in L^2(\R^N ), \, f_1(x,u_1), f_2(x,u_2) \in L^2(\\R^N), \leqno(i)$$
$$ \exists \, \kappa>0  \; s.t.  \, f_1(x,u_1), f_2(x,u_2)  \, \geq \, \kappa \phi(x)  \; \forall  u_1,u_2 \; \in\,  L^2(\R^N ), \leqno(ii)$$
$$ \exists K > \kappa>0 \; s.t.  \, f_1(x,u_1), f_2(x,u_2) \, \leq \, K \phi(x)  \; \forall  u_1, u_2 \, \in\,  L^2(\R^N ). \leqno(iii)$$
  $(H'_{F}):$ $\;$ 
$\frac{f_1(x,u_1)}{|u_1|}$ and $\frac{f_2(x,u_2)}{|u_2|}$ are decreasing w.r.t. $u_1$ and $u_2$. 

\vskip0.2cm \noindent
  We introduce 2 sets : 
\begin{align*}{\cal K_S}^+: = \{(u_1,u_2)\in X^2\, / \,& u_1 \in \big(\frac{\kappa y_1\phi}{\max(y_1,y_2)(\Lambda^* - \mu)}\, , \, \frac{Ky_1\phi}{\min(y_1,y_2)(\Lambda^* - \mu)} \big) ,\\
 & u_2 \in \big(\frac{\kappa y_2\phi}{\max(y_1,y_2)(\Lambda^* - \mu)}\, , \, \frac{Ky_2\phi}{\min(y_1,y_2)(\Lambda^* - \mu)} \big)\}
\end{align*}
 for $ \mu < \Lambda^*  $, and 
\begin{align*}{\cal K_S}^-: = \{(u_1,u_2)\in X^2\, / \,& u_1 \in \big(\frac{K y_1\phi}{\min(y_1,y_2)(\Lambda^* - \mu)}\, , \, \frac{\kappa y_1\phi}{\max(y_1,y_2)(\Lambda^* - \mu)} \big) ,\\
 & u_2 \in \big(\frac{K y_2\phi}{\min(y_1,y_2)(\Lambda^* - \mu)}\, , \, \frac{\kappa y_2\phi}{\max(y_1,y_2)(\Lambda^* - \mu)} \big)\}
\end{align*}
 for $ \Lambda^* < \mu $.

\begin{thm} \label{T4}
 If $(H_A)$ and  $(H_{F})$ are  satisfied there exists $\delta >0$,  depending on $f_1$ and $f_2$ such that if 
$\Lambda^*  - \delta <\mu <\Lambda^* $  (resp.
 $\Lambda^* < \mu < \Lambda^*  + \delta$),  (with  $  \delta <  min\{ \frac{\xi_2 - \xi_1}{2},\lambda_2 -\Lambda\} $)    
 System $(S)$ has a solution which is in $\cal K_S^+$, (resp. in $\cal K_S^-$).
 Moreover, if $(H'_{F})$ is satisfied, the solution is unique. 
\end{thm}
{\bf Proof of Theorem \ref{T4}:}$\;$ We use of course the results above as well as previous results for linear systems obtained in \cite{AF2017} where
 Theorem \ref{T4} is shown for suitable assumptions on $f_1$ and $f_2$ ( independent  on $u$).  
\vskip0.2cm  \noindent
{\bf 1. Maximun and anti-maximum principles}
$\;$
\par \noindent
We  diagonalize System$(S)$  thanks to the change of basis  matrix $P$,  and we get  a  system of 2 equations. 
Here
$$P =\left(\begin{array}{cccc}
b&b\\\xi_1 - a &\xi_2 - a
 \end{array}
 \right)\, , \;   P^{-1}  =\frac{1}{ b(\xi_1 - \xi_2)} \left(\begin{array}{cccc}
 a-\xi_2  &  b\\  \xi_1 -a & -b
 \end{array}
 \right)\, , \; $$ 
Set 
\begin{equation} \label{matrix} D:=P^{-1} A P = \left(\begin{array}{cccc}
\xi_1 &0 \\ 0& \xi_2
 \end{array}
 \right)\,; \; U=PV\,; \; G:=  P^{-1}F.\end{equation} 
 We obtain
\begin{equation}\label{S'} LV = DV+\mu V + G \end{equation}
which is a system of 2 equations (with obvious notation):
$$Lv_1 =(\xi_1 + \mu) v_1 + g_1(u_1,u_2);$$
$$Lv_2=(\xi_2 + \mu) v_2 + g_2(u_1,u_2).$$
Note that $g_1$ and $g_2$ are in $X$.
\par \noindent
The second equation, where the parameter $\xi_2+ \mu$ stays away (below)  from $\Lambda$, has a $\phi$ bounded solution $v_2$. 
Concerning the first equation, we apply Theorem \ref{T1} above.  We compute $g_1$, $g_2$  and get 
\par \noindent
$$ \exists \, \kappa'>0  \; s.t.  \, g_1(x,u_1,u_2) \, \geq \, \kappa' \phi(x)  \; \forall  u_1,u_2 \; \in\,  L^2(\R^N ), \leqno(ii')$$
$$ \exists K' > \kappa'>0 \; s.t.  \, g_1(x,u_1,u_2), \, |g_2(x,u_1,u_2)| \, \leq \, K' \phi(x)  \; \forall  u_1, u_2 \, \in\,  L^2(\R^N ), \leqno(iii')$$
 where $\kappa'$ and $K'$ are 2 positive constants depending on $\kappa$, $K$ and 
on the coefficients of $A$. 
This follows from  $\xi_1-\xi_2>0$ and $(a-\xi_2) = \frac{a-d}{2} + \frac{\sqrt{(a-d)^2+4bc}}{2}$ with 
$(a-d)^2+4bc > (a-d)^2$, so that  
$$g_1=\frac{1}{\xi_1-\xi_2} [(a-\xi_2) f_1 + bf_2]> \kappa' \phi >0. $$
Analoguously we have $g_1<K'\phi$. Therefore Theorem 2 holds here with $\delta = \min(\delta_0, \frac{ \kappa'}{c_0 K'},\frac{\xi_1 - \xi_2}{2} )$.
Finally we deduce from the maximum principle for $\Lambda^*-\delta <\mu < \Lambda^*$ that $v_1> \frac{\kappa'}{\Lambda^* - \mu} \phi>0$. 
\par \noindent
If $ \Lambda^*<\mu < \Lambda^*+\delta$, reasoning similarly, we deduce  $v_1<  \frac{K'}{\Lambda^* - \mu} \phi<0$.
As $\mu \rightarrow \Lambda^*$, $v_1$ tends to $\infty$ when $v_2$ stays bounded. Indeed, by Remark \ref{RT1}, 
$$\|v_2\|_X \leq \frac{K'}{|\Lambda - \xi_2 - \mu|} + 2 c_0 K'  < \frac{2K'}{\xi_1 - \xi_2} +2c_0 K';$$
the last inequality follows from $\delta < \frac{\xi_1 - \xi_2}{2}$.
\par \noindent
Now we go back to $U=PV$. 
$$u_1=av_1 + bv_2\, , \; u_2=(\xi_1 - a) v_1 + (\xi_2 -a)v_2.$$
Combining the estimates above on $v_1$ and $v_2$, we conclude that, as $|\Lambda^* - \mu| \rightarrow 0$, there exists 
$\delta^*$, depending only on $L, A, \kappa, K$ such that 
 as  $\mu \nearrow \Lambda^*$, $u_1$ has the sign of $a >0$ and $u_2>0$. 
If  $\mu \searrow \Lambda^*$, $u_1$ has the sign of $-a<0$ and $u_2<0$.

\vskip0.2cm  \noindent
{\bf 2. Existence of the solution in $\cal K_S^+$, (resp. in $\cal K_S^-$ )}
$\;$
\par \noindent
\underline{Sub-supersolutions: } $\;$ 
\par \noindent
$1$. Case $\Lambda^* - \delta^*<\mu<\Lambda^*$. Recall  that $Y$  has positive  components  $y_1$ and $y_2$, and 
the principal eigenvector   $\Phi =Y \phi$  satisfies  
$$ L\Phi -\mu \Phi - A\Phi= (\Lambda^* - \mu )\Phi.$$
 Inspired by the case of one equation, we seek a subsolution $U_0$ of the form $cY\Phi$. 
$$L(U-U_0) =   A(U-U_0) + \mu (U-U_0) +(F(x,U)-(\Lambda^* - \mu)c\Phi).$$
For $c$ such that $F(x,U) -(\Lambda^* - \mu)cY \phi(x)>0$,  for $\mu < \Lambda^*$,  we get $U-U_0 >0$ by the maximum principle. 
Finally, since $ F(x,U) - \dfrac{\kappa}{\max(y_1,y_2)}  Y\phi>0  $,  a subsolution is 
$$U_0 = \frac{\kappa}{\max(y_1,y_2)}  \frac{1}{(\Lambda^* - \mu)} Y\phi .$$

\vskip0.2cm \noindent
Analogously $U^0 = \dfrac{K}{{\min(y_1,y_2)} (\Lambda^*- \mu)}Y \phi$ is a supersolution.
\vskip0.3cm \noindent
$2$. Case $\,$ $\Lambda^*<\mu<\Lambda^*+\delta^*$.
We have similar results with change of sign and replacing $K$ by $\kappa$. 
$$V_0 = \frac{K}{{\min(y_1,y_2)} (\Lambda^*- \mu)}Y\phi$$
$$V^0 = \frac{\kappa}{\max(y_1,y_2)}  \frac{1}{(\Lambda^* - \mu)} Y\phi $$

\vskip0.2cm  \noindent
\underline{ The operator $T$: }$\;$
We define $T: ( u_1,u_2) \longrightarrow( w_1,w_2)$ where $(w_1,w_2)$ is the solution to the linear system
$$ \left \{ \begin{array}{lcl}
   Lw_1 \,&=&\, (a+\mu) w_1 + bw_2 +f_1(x,u_1)\\
 Lw_2 \,&=&\, cw_1+(d+\mu) w_2 +f_2(x,u_2)        
          \end{array} \right . \, \; in \; \R^N ,.\leqno(S') $$
$$ \; w_1(x), w_2(x)_{|x|\rightarrow \infty} \rightarrow 0 .$$
\underline{The rectangle:}$\;$
If $(u_1,u_2) \in {\cal K_S^+}$ for $\Lambda^* - \delta^* <\mu<\Lambda^*$ (resp. $(u_1,u_2) \in {\cal K_S^-}$ for $\Lambda^* < \mu < \Lambda^* + \delta^* $)
 then $(w_1,w_2)\in {\cal K_S^+}$ (resp ${\cal K_S^-}$). 
Indeed, for $\Lambda^* - \delta^*\mu<\Lambda^*$, this can be written with obvious notations 
$$L(W- U_0)=(\mu+A)(W-U_0)+F;$$ for $\mu < \Lambda^*$, 
since $F$ has non negative components, $F\not \equiv 0$, then $W-U_0>0$. 
Analogously, we obtain the supersolution $U^0-W>0$. 

We argue exactly as for one equation: $\cal K_S^+$ or $\cal K_S^-$ is invariant by $T$ and 
$LW=(A+\mu)W+F(x,U)$ can be written 
 $ W= (L-A-\mu I)^{-1}\hat F(x,u) = T(U)$. Since by \cite{ATa2007} ,\cite{AFTa2007},  the resolvent 
$ R(\mu):= (L-\mu I)^{-1}$ is compact  in $X$ for $\mu \in (\Lambda^* - \delta^*, \Lambda^*)$ or 
$(\Lambda^*, \Lambda^* + \delta^*)$,   and  since $\hat F:u \rightarrow F(x,u)$ is continuous, $T= R(\mu)\hat F$ is compact.

 We apply the fixed point theorem. There exists a solution $U$.
\vskip0.2cm  \noindent
{\bf 3. Uniqueness}
$\;$
\par \noindent

We assume now $(H'_{F})$. 
assume there are 2 positive  solutions $(u_1,u_2)$ and $(v_1,v_2)$ to $(S)$; 
for the first equation we have
$Lu_1=(\mu+a)u_1+bu_2+f_1(x,u_1)$  and $Lv_1=(\mu+a)v_1+bv_2+f_1(x,v_2)$.
Since we are in $\cal K^+$ (resp. $\cal K ^-$), divide by $bu_1$ the first equation and by $bv_1$ the second one and subtract:
\begin{equation} \label{E1} \frac{-\Delta u_1}{bu_1} +\frac{\Delta v_1}{bv_1} =
\frac{u_2}{u_1} - \frac{v_2}{v_1} + \frac{f_1(x,u_1)}{bu_1} -  \frac{f_1(x,v_1)}{bv_1}.\end{equation}
Exactly as in \cite{BrOs} multiply by $(u_1^2-v_1^2)$ and integrate; 
hence 
\[ \int(\frac{-\Delta u_1}{bu_1} +\frac{\Delta v_1}{bv_1}) (u_1^2-v_1^2)=
\int (\frac{u_2}{u_1} - \frac{v_2}{v_1} + \frac{f_1(x,u_1)}{bu_1} -  \frac{f_1(x,v_1)}{bv_1})(u_1^2-v_1^2).\]
The first terme is non-negative by  ( \ref{brezis}): 
$$\int (\frac{-\Delta u_1}{bu_1} +\frac{\Delta v_1}{bv_1})(u_1^2-v_1^2)>0.$$
 Then do exactly the same calculus with the second equation in $(S)$  and add these two lines: we derive from (\ref{E1}) that  $T_1=T_2$  with
\begin{align*}
T_1 & =  \int (\frac{-\Delta u_1}{bu_1} +\frac{\Delta v_1}{bv_1})(u_1^2-v_1^2) + \int (\frac{-\Delta u_2}{cu_2} +\frac{\Delta v_2}{cv_2})(u_2^2-v_2^2).\\
T_2 & =\int (\frac{u_2}{u_1} - \frac{v_2}{v_1} + \frac{f_1(x,u_1)}{bu_1} -  \frac{f_1(x,v_1)}{bv_1})(u_1^2-v_1^2) +\\
&\quad \quad \quad \quad \quad \quad \quad \quad \quad \quad \quad \quad \quad \quad \quad  \int (\frac{u_1}{u_2} - \frac{v_1}{v_2} +
 \frac{f_2(x,u_2)}{cu_2} -  \frac{f_2(x,v_2)}{cv_2})(u_2^2-v_2^2).
\end{align*}
Of course the 1st term $T_1$  is non-negative by (\ref{brezis}). 
By $(H'_{F})$, $$\int ( \frac{f(x,u_1)}{bu_1} -  \frac{f_1(x,v_1)}{bv_1})(u_1^2-v_1^2) + \int
 (\frac{f_2(x,u_2)}{cu_1} -  \frac{f_2(x,v_2)}{cv_1})(u_2^2-v_2^2) <0. $$
We develop what is left and get 
\begin{align*}
\int (\frac{u_2}{u_1} - \frac{v_2}{v_1})(u_1^2-v_1^2)+  \int &(\frac{u_1}{u_2} - \frac{v_1}{v_2} )
(u_2^2-v_2^2)
= \\
&- \int (\sqrt{\frac{u_2v_1^2}{u_1}} - \sqrt{\frac{u_1v_2^2}{u_2}})^2- \int (\sqrt{\frac{v_2u_1^2}{v_1}} - \sqrt{\frac{v_1u_2^2}{v_2}})^2<0
\end{align*}
Hence $T_1 = T_2 = 0$ and $u_1=v_1$,$ u_2=v_2$. The solution is unique.

{}

Authors: \\
B.Alziary  \\
Institut de Math\'ematiques de Toulouse CeReMath, TSE, Univ. Toulouse 1, \\
J.Fleckinger \\
Institut de Math\'ematiques de Toulouse CeReMath, Univ. Toulouse 1, \\
address: 21 all\'ee de Brienne 31042 Toulouse Cedex \\
email: $\;$  alziary@ut-capitole.fr,$\,$  jfleckinger@gmail.com

\end{document}